\documentclass[a4paper,10pt]{article}

\usepackage{setspace}
\usepackage[bottom]{footmisc}
\usepackage{color}
\usepackage{subfig}
\usepackage{hyperref}
\usepackage{multirow}
\usepackage{cite}
\usepackage{amsmath} 
\usepackage{mathtools}
\usepackage{cleveref}
\usepackage{graphicx}
\usepackage{epsfig}
\usepackage{epstopdf}
\usepackage{float}
\usepackage{titling}
\usepackage[margin=1.0in]{geometry}
\usepackage{etoolbox}
 
\title{Study of Mach reflection in inviscid flows}
\makeatletter
\renewcommand\@date{{%
  \vspace{-\baselineskip}%
  \large\centering
  \begin{tabular}{@{}c@{}}
    S R Siva Prasad Kochi\textsuperscript{1} \\
    \normalsize siva.ksr@gmail.com
  \end{tabular}%
  \quad and\quad
  \begin{tabular}{@{}c@{}}
    M Ramakrishna\textsuperscript{2} \\
    \normalsize krishna@ae.iitm.ac.in
  \end{tabular}

  \bigskip

  \textsuperscript{1}Dept. of Aerospace Engg., IIT Madras.\par
  \textsuperscript{2}Professor, Dept. of Aerospace Engg., IIT Madras.

  \bigskip

  \today
}}
\makeatother
\begin{document}

\maketitle

\begin{abstract}
 In this paper, we study the Mach reflection phenomenon in inviscid flows using a higher order discontinuous Galerkin method and overset grids. We use the shock capturing procedure proposed in \cite{srspkmr1} using overset grids to capture the discontinuities occurring in the supersonic flow over a wedge accurately. In this procedure, we obtain a coarse grid solution first and using the troubled cell data, we construct an overset grid which is approximately aligned to all the discontinuities. We rerun the solver with the coarse grid solution as the initial condition while using the troubled cell indicator and the limiter only on the overset grid. This allows us to capture the discontinuities accurately. Using this procedure, we have obtained the solution for Mach $3.0$ and $4.0$ flow over a wedge for various wedge angles and determined the detachment criterion and the Von Neumann condition accurately. We have also determined the Mach stem height for various wedge angles for these Mach numbers. We have also demonstrated the hysteresis that occurs in the transition from regular reflection to Mach reflection.
 
 {{\bf Keywords:} Mach reflection, discontinuous Galerkin method, overset grids, hysteresis \\}
\end{abstract}

\section{Introduction}\label{sec:intro}

\begin{figure}[h!]
\begin{center}
\scalebox{0.85}{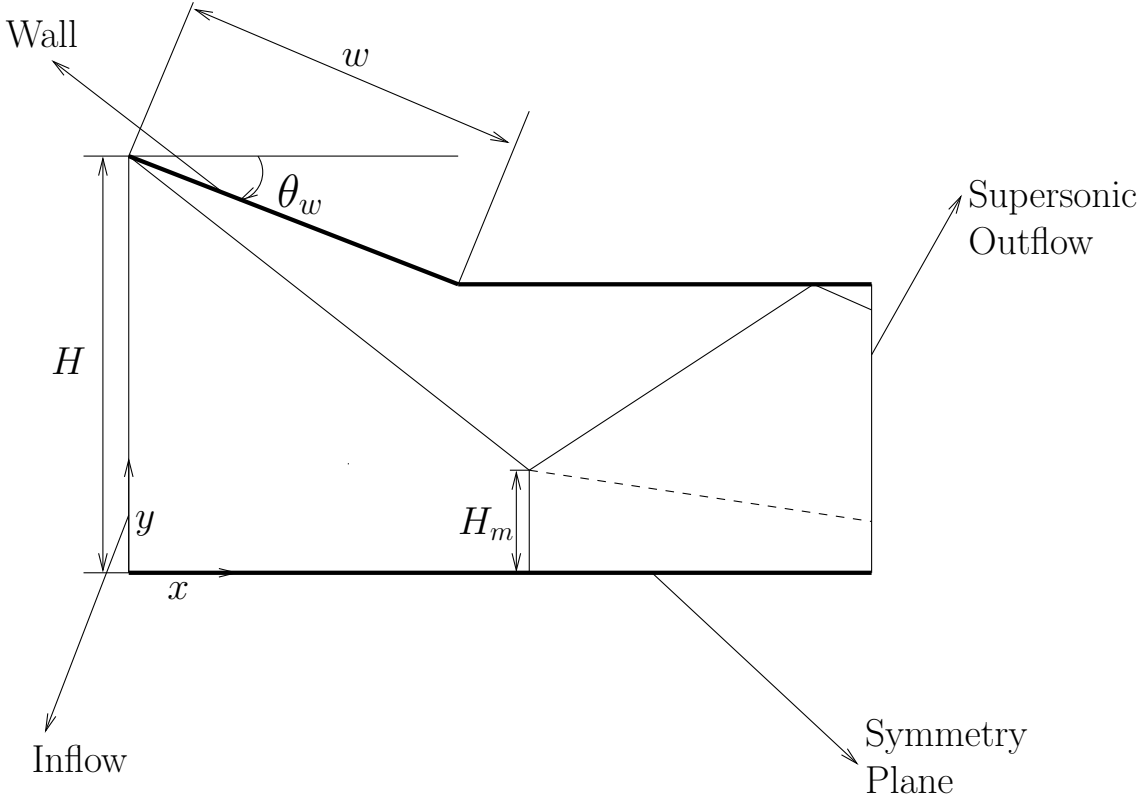}
\caption{Computational domain for Mach reflection showing the shock structure and the Mach stem height ($H_{m}$)}
\label{fig:CompDomain}
\end{center}
\end{figure}

\noindent In this paper, we study the transition between regular reflection (RR) and Mach reflection (MR) of steady shock waves in inviscid flows using discontinuous Galerkin method (DGM) along with overset grids. This has been studied quite extensively in literature experimentally \cite{cplb}, numerically \cite{vzb}, and analytically \cite{bendor}. We use the shock capturing procedure proposed in \cite{srspkmr1} using overset grids and a higher order method (DGM) to capture the discontinuities occurring in the supersonic flow over a wedge accurately. In this procedure, we obtain a coarse grid solution first and using the solution and the troubled cell data, we construct an overset grid which is approximately aligned to all the discontinuities. We rerun the solver with the coarse grid solution as the initial condition while using the troubled cell indicator and the limiter only on the overset grid. This allows us to capture the discontinuities accurately. We solve the Euler equations in the computational domain shown in Figure \ref{fig:CompDomain} using DGM and overset grids for Mach 3.0 and 4.0 flow over a wedge for different wedge angles to determine the transition between RR and MR. In this way, we determine the detachment criterion and the Von Neumann condition accurately to demonstrate the hysteresis that occurs in such a flow. 
\\

\noindent The paper is organized as follows. We describe the formulation of the discontinuous Galerkin method used for all our results in Section \ref{sec:dgm}, the procedure used for shock capturing using overset grids is described in Section \ref{sec:overset}, the results are described in Section \ref{sec:Results} and we conclude the paper in Section \ref{sec:conc}.

\section{Description of discontinuous Galerkin method}\label{sec:dgm}

\noindent Consider the Euler equations in conservative form as given by

\begin{equation}\label{2dEulerEquations}
\frac{\partial \textbf{Q}}{\partial t} + \frac{\partial \textbf{F(Q)}}{\partial x} + \frac{\partial \textbf{G(Q)}}{\partial y} = 0 \quad \text{in the domain} \quad \Omega
\end{equation}
\noindent where $\textbf{Q} = (\rho, \rho u, \rho v, E)^{T}$, $\textbf{F(Q)}=u\textbf{Q} + (0, p, 0, pu)^{T}$ and $\textbf{G(Q)}=v\textbf{Q} + (0, 0, p, pv)^{T}$ with $p = (\gamma -1)(E-\frac{1}{2}\rho (u^{2}+v^{2}))$ and $\gamma = 1.4$. Here, $\rho$ is the density, $(u,v)$ is the velocity, $E$ is the total energy and $p$ is the pressure. We approximate the domain $\Omega$ by $K$ non overlapping elements given by $\Omega_{k}$. 
\\
\\
\noindent We look at solving \eqref{2dEulerEquations} using the discontinuous Galerkin method. We approximate the local solution in an element $\Omega_{k}$, where $k$ is the element number, as a polynomial of order $N$ which is given by:

\begin{equation}\label{modalApprox}
 Q_{h}^{k}(r,s) = \sum_{i=0}^{N_{p}-1} Q_{i}^{k} \psi_{i}(r,s)
\end{equation}

\noindent where $N_{p}=(N+1)(N+1)$ and $r$ and $s$ are the local coordinates. Here, the subscript $i$ represents the particular degree of freedom, $h$ represents the grid size, and the superscript $k$ is the element number. The polynomial basis used ($\psi_{i}(r,s)$) is the tensor product orthonormalized Legendre polynomials of degree $N$. The number of degrees of freedom are given by $N_{p}=(N+1)(N+1)$. Now, using $\psi_{j}(r,s)$ as the test function, the weak form of the equation \eqref{2dEulerEquations} is obtained as

\begin{equation}\label{weakFormScheme}
 \sum_{i=0}^{N_{p}-1} \frac{\partial Q_{i}^{k}}{\partial t} \int_{\Omega_{k}} \psi_{i} \psi_{j} d\Omega + \int_{\partial \Omega_{k}} \hat{F} \psi_{j} ds - \int_{\Omega_{k}} \vec{F} \cdot \nabla \psi_{j} d\Omega = 0 \quad j = 0,\ldots,N_{p}-1
\end{equation}

\noindent where $\partial \Omega_{k}$ is the boundary of $\Omega_{k}$, $\vec{F} = (\textbf{F(Q)},\textbf{G(Q)})$ and $\hat{F} = \bar{F^{*}}\cdot\hat{n}$ where $\bar{F^{*}}$ is the monotone numerical flux at the interface which is calculated using an exact or approximate Riemann solver and $\hat{n}$ is the unit outward normal. This is termed to be $\mathbf{P}^{N}$ based discontinuous Galerkin method.
\\
\\
\noindent Equation \eqref{weakFormScheme} is integrated using an appropriate Gauss Legendre quadrature and is discretized in time by using the fifth order Runge-Kutta time discretization given in \cite{butcher} unless otherwise specified. To control spurious oscillations which occur near discontinuities, a limiter is used with a troubled cell indicator. We have used the KXRCF troubled cell indicator \cite{kxrcf} and the compact subcell WENO (CSWENO) limiter proposed in \cite{srspkmr2} for all our calculations.

\section{Overset grids and shock capturing}\label{sec:overset}

\begin{figure}[htbp]
\begin{center}
\includegraphics[scale=1.35]{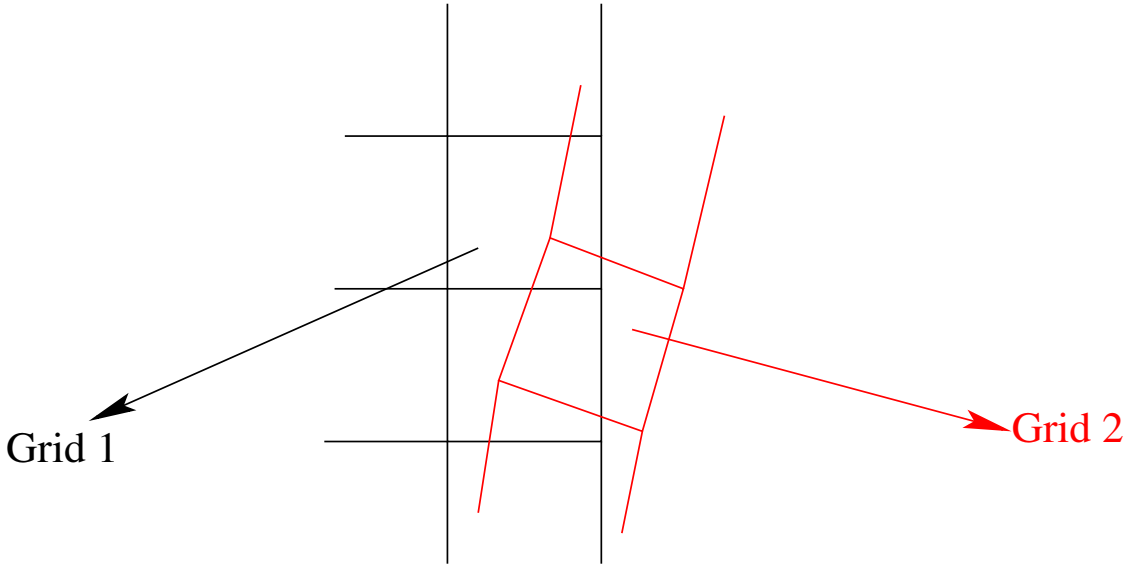}
\caption{Two overlapping grids (Grid 1 in black and Grid 2 in Red)}
\label{fig:OversetGridExpl}
\end{center}
\end{figure}

\noindent Overset grids consist of multiple grids which overlap each other as shown in Figure \ref{fig:OversetGridExpl}. When using DGM on overset grids, there are two possible approaches to handle data communication between the grids. One is a face based communication approach developed in \cite{gbot}, where solutions at an overset interface are obtained from the donor element, and then the boundary condition is applied weakly by imposing a numerical flux at the flux interpolation points. The other is an element based communication approach developed in \cite{nsm}, where the internal degrees of freedom of cells near the overset interface are obtained from the donor element. We use the new element based communication approach developed in \cite{srspkmr3} for the data communication between the grids.

\begin{figure}[htbp]
\begin{center}
\includegraphics[scale=0.22]{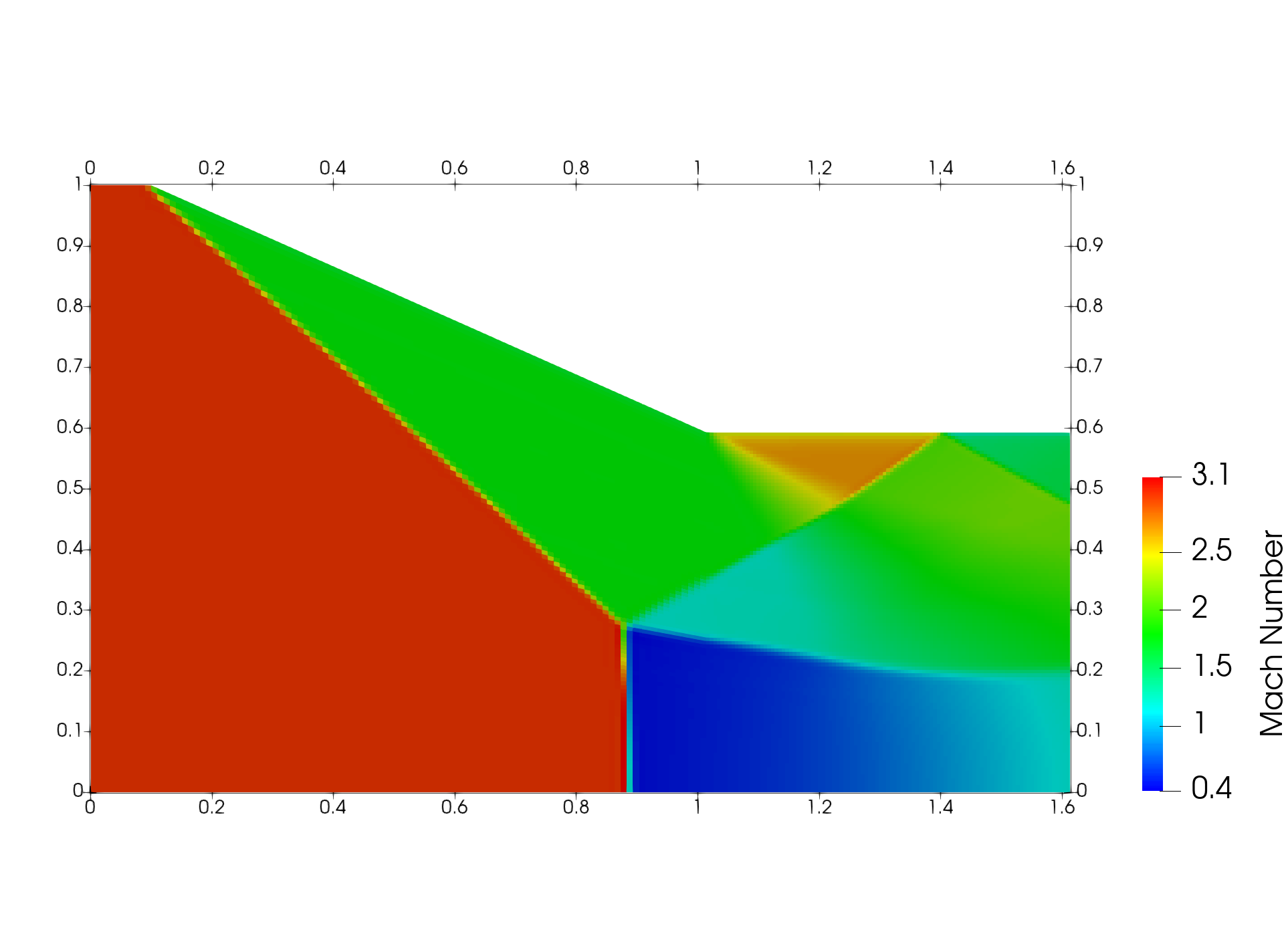}
\caption{Coarse grid solution for Mach 3.0 flow over a $24^{\circ}$ wedge using $\mathbf{P}^{1}$ based discontinuous Galerkin method}
\label{fig:M3p0Angle24DCoarseSol}
\end{center}
\end{figure}

\begin{figure}[htbp]
\begin{center}
\includegraphics[scale=0.22]{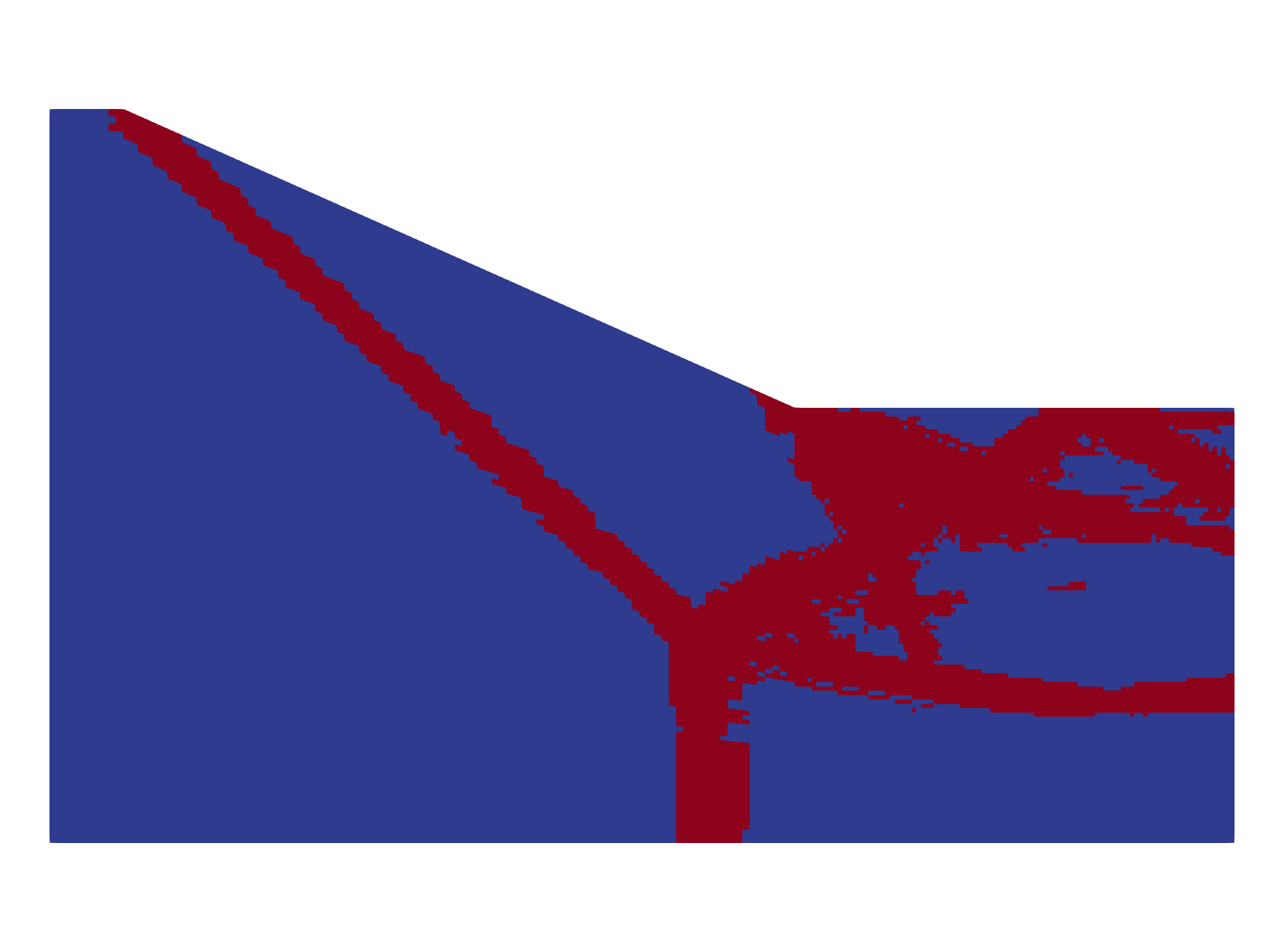}
\caption{Troubled cells obtained using the KXRCF troubled cell indicator \cite{kxrcf} for Mach 3.0 flow over a $24^{\circ}$ wedge}
\label{fig:M3p0Angle24DTroubledCells}
\end{center}
\end{figure}

\begin{figure}[htbp]
\begin{center}
\includegraphics[scale=0.22]{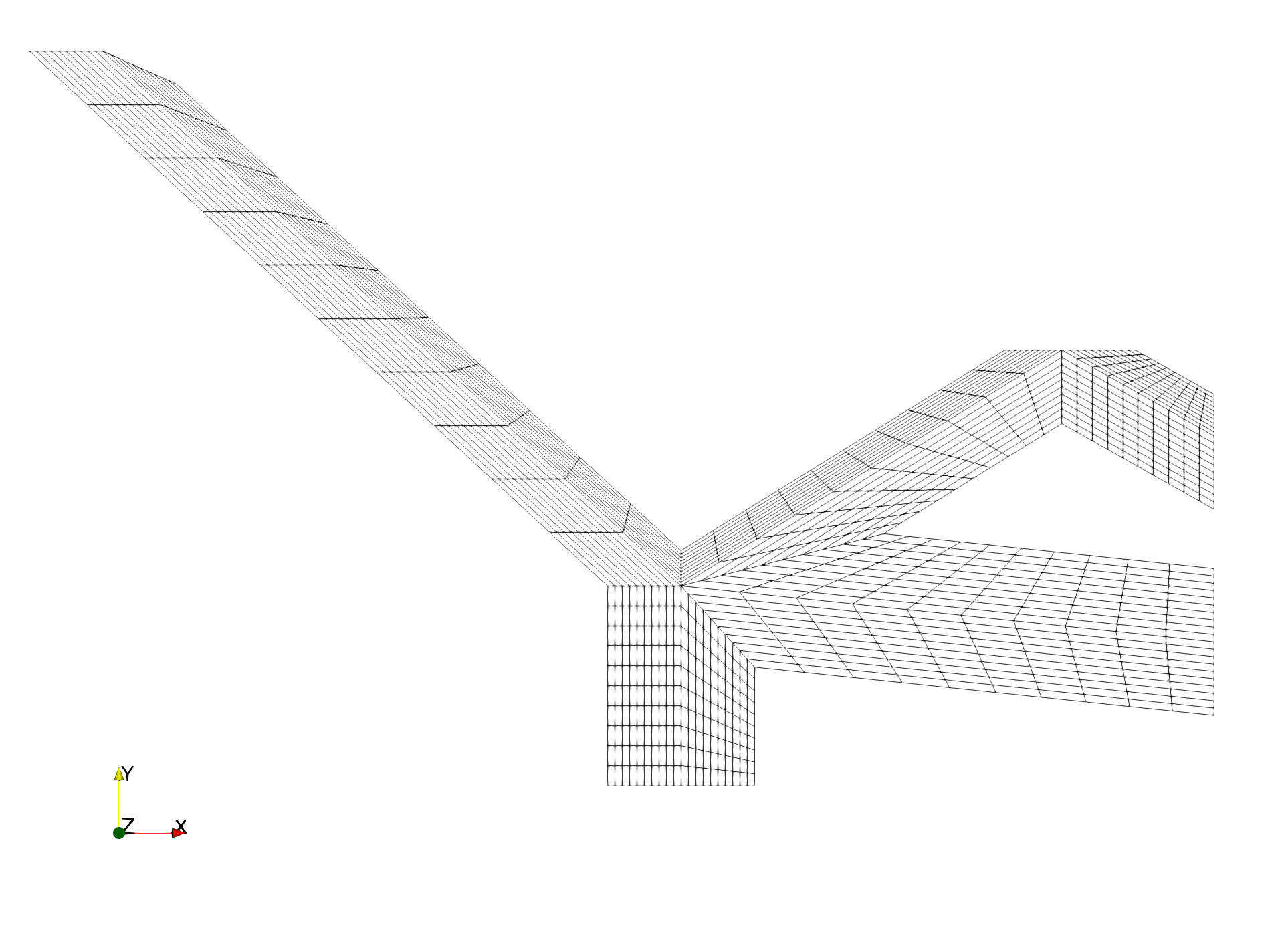}
\caption{Example overset grid for Mach 3.0 flow over a $24^{\circ}$ wedge constructed such that the discontinuities are approximately parallel to a grid line}
\label{fig:M3p0Angle24DOversetGrid}
\end{center}
\end{figure}

\noindent For capturing the shocks accurately, we use the procedure developed in \cite{srspkmr1} using overset grids. A brief explanation of the procedure is given below:
\\
\\
\noindent \textbf{Step 1:} Run the solver on a coarse grid with a given troubled cell indicator and limiter to steady state and obtain the solution. As an example, we show the coarse grid solution obtained for Mach 3.0 flow over a $24^{\circ}$ wedge in Figure \ref{fig:M3p0Angle24DCoarseSol}
\\
\\
\noindent \textbf{Step 2:} Look at the troubled cells to locate the discontinuities (shocks) that occur in the solution. The troubled cell profile obtained for the Mach 3.0 flow over a $24^{\circ}$ wedge using the KXRCF troubled cell indicator \cite{kxrcf} is shown in Figure \ref{fig:M3p0Angle24DTroubledCells}. From this figure, we can see that the troubled cells give us a good idea of the location of the shocks, the contact discontinuity and the initial expansion that forms near the wedge.
\\
\\
\noindent \textbf{Step 3:} Construct an overset grid conforming to the computational domain which is refined in a direction perpendicular to the discontinuities such that they are approximately parallel to a grid line. This overset grid also encompasses all the troubled cells. An example overset grid constructed in such a fashion for the Mach 3.0 flow over a $24^{\circ}$ wedge is shown in Figure \ref{fig:M3p0Angle24DOversetGrid}.
\\
\\

\begin{figure}[htbp]
\begin{center}
\includegraphics[scale=0.22]{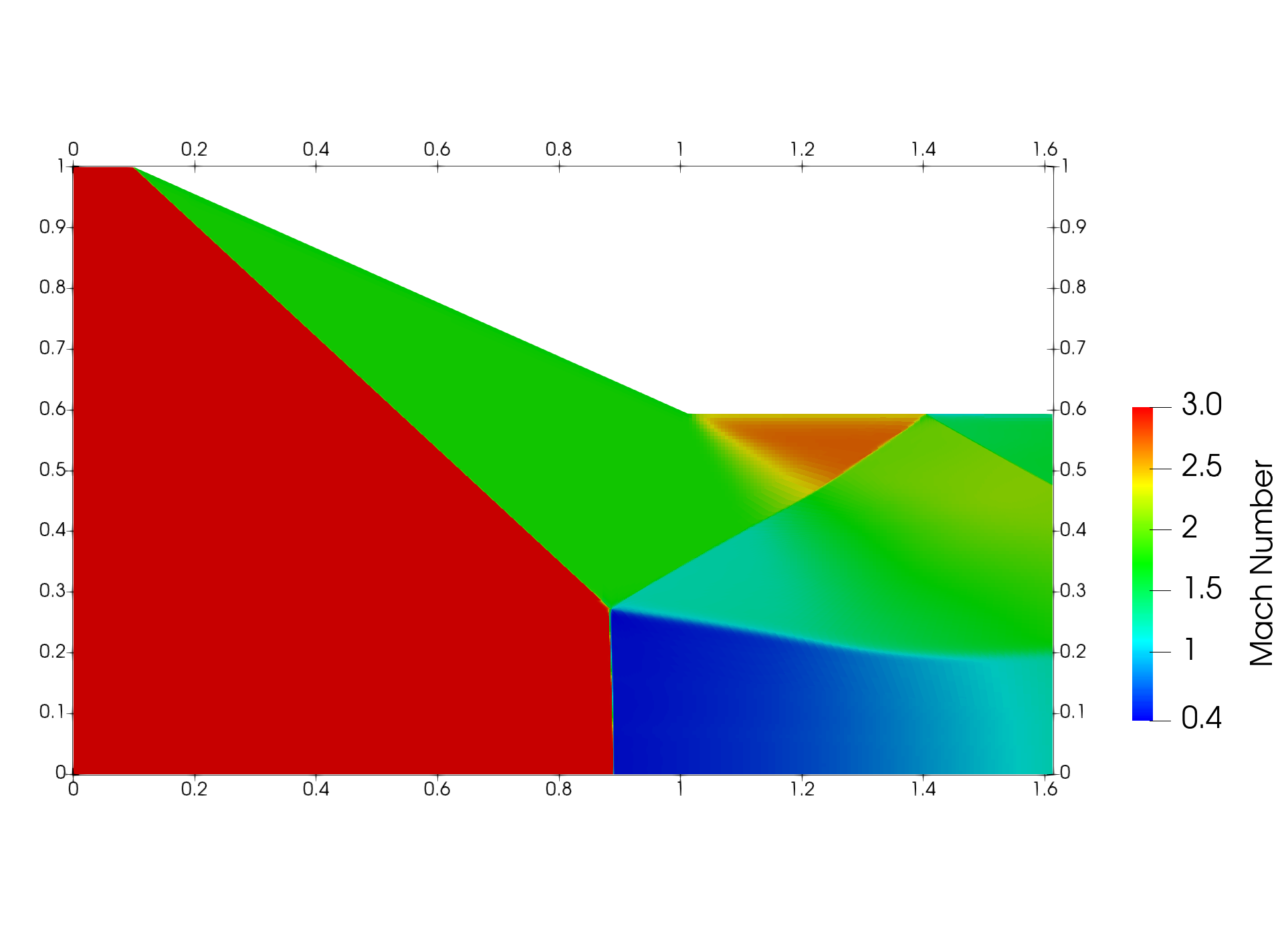}
\caption{Final solution obtained for Mach 3.0 flow over a $24^{\circ}$ wedge using an overset grid and $\mathbf{P}^{4}$ based discontinuous Galerkin method with discontinuities captured accurately}
\label{fig:M3p0Angle24DFinal}
\end{center}
\end{figure}

\noindent \textbf{Step 4:} Using this overset grid, we rerun the solver with the coarse grid solution as the initial condition. While running the solver, we use the troubled cell indicator and the limiter only on the overset grid. We also use a high resolution numerical flux on the overset grid to capture the shock accurately. We have used the SLAU2 \cite{ks3} numerical flux in the overset grid and the less expensive Lax-Friedrichs flux elsewhere. Using this procedure, we obtain a more accurate solution with the discontinuities approximately aligned to a grid line. The final solution obtained in this fashion for Mach 3.0 flow over a $24^{\circ}$ wedge is shown in Figure \ref{fig:M3p0Angle24DFinal}.
\\
\\

\noindent Regarding the computational cost of the scheme proposed for the current problem, we note that the troubled cell data shown in Figure \ref{fig:M3p0Angle24DTroubledCells} are obtained from a coarse grid of 20000 elements using $\mathbf{P}^{1}$ based DGM after converging to a residue of about $1\times 10^{-6}$ which happens in about 8000 iterations. This coarse grid solution is not a fully converged solution but is good enough to be used as an initial condition for the overset grid higher order ($\mathbf{P}^{4}$ based DGM) solution. The time taken per iteration per degree of freedom using $\mathbf{P}^{1}$ based DGM for this problem is about $1.084\times 10^{-6}$s. This tells us that the coarse grid troubled cell data is obtained in about $693.76$s which is a little less than 12 minutes. Most of the work involved is in constructing the overset grid which is done manually. After the overset grid is constructed and the solver is run on the new grid with the coarse grid solution as the initial condition, the residue converges to $1\times 10^{-16}$ in about 32000 iterations. For $\mathbf{P}^{4}$ based DGM, the time taken per iteration per degree of freedom is about $1.345 \times 10^{-6}$s for this problem. This tells us that the final solution (which is fifth order accurate) is obtained in about $21520$s which is about six hours. All these calculations are done on a 3.60 GHz, Intel(R) Core(TM) i7-7700 CPU with a single thread. The calculations for remaining angles are quite similar. We also note that we have a parallelised code and we obtain the solution much faster based on the number of threads used.

\section{Results}\label{sec:Results}

\begin{figure}[htbp]
\begin{center}
\includegraphics[scale=0.22]{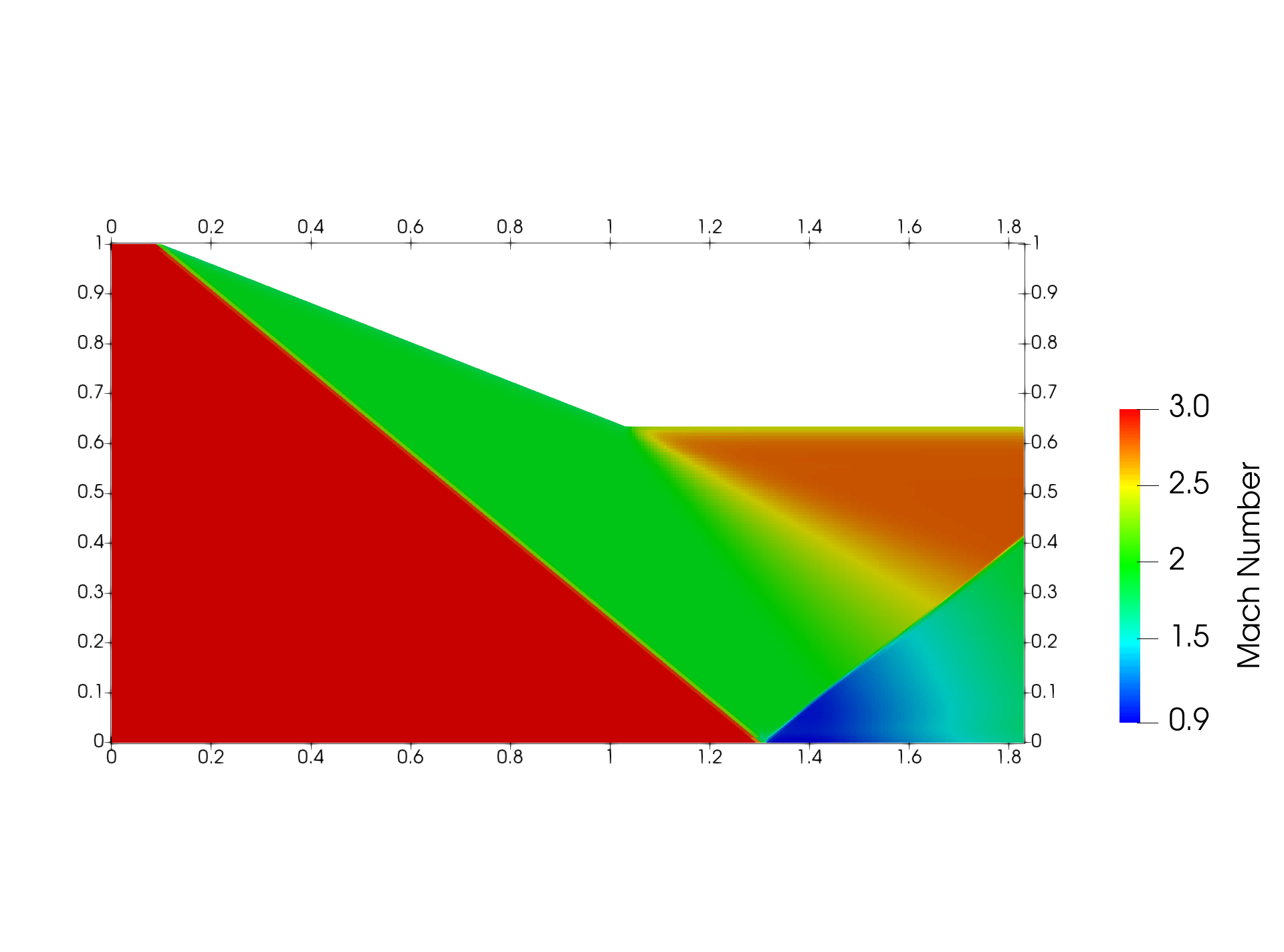}
\caption{Final solution obtained for Mach 3.0 flow over a $21.45^{\circ}$ wedge with an overset grid and $\mathbf{P}^{4}$ based discontinuous Galerkin method with discontinuities captured accurately using an impulsive start}
\label{fig:M3p0Angle21p45DFinal}
\end{center}
\end{figure}

\begin{figure}[htbp]
\begin{center}
\includegraphics[scale=0.22]{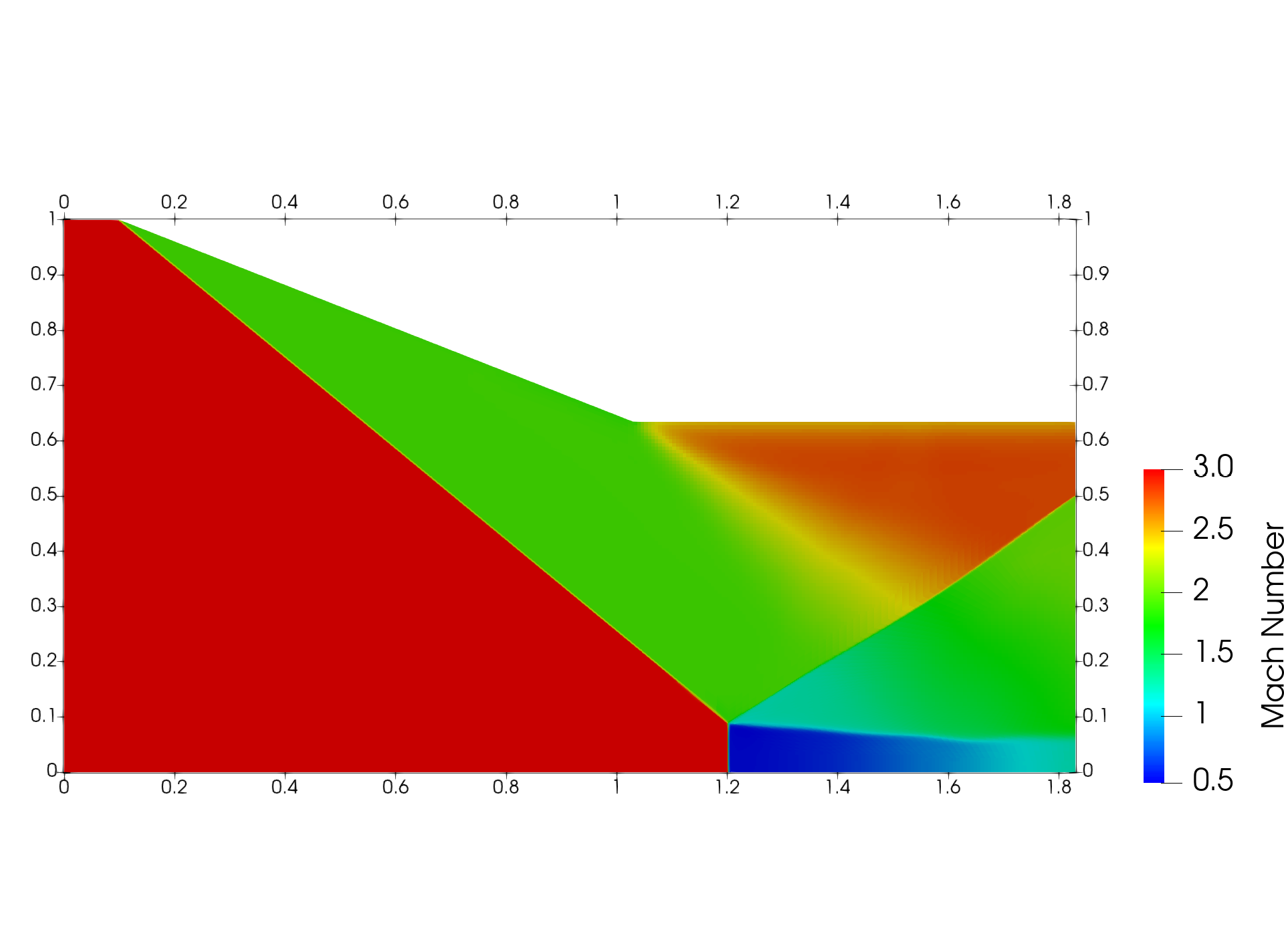}
\caption{Final solution obtained for Mach 3.0 flow over a $21.46^{\circ}$ wedge with an overset grid and $\mathbf{P}^{4}$ based discontinuous Galerkin method with discontinuities captured accurately using an impulsive start}
\label{fig:M3p0Angle21p46DFinal}
\end{center}
\end{figure}

\textbf{1) Mach number 3.0:} We solve the two-dimensional Euler equations given by \eqref{2dEulerEquations} using the discontinuous Galerkin method for various wedge angles between $\theta_{w}=19.5^{\circ}$ and $\theta_{w}=24^{\circ}$ near the transition criterion (transition from Regular reflection to Mach reflection) for Mach number 3.0 and $w/H=1.0$ in the computational domain shown in Figure \ref{fig:CompDomain}. We consider two cases to demonstrate the hysteresis phenomenon. We solve the equations using an impulsive start as the initial condition for the first case, and the converged solution for $\theta_{w}=24^{\circ}$ as the second case. The first case is a numerical model to obtain the detachment criterion and the second case demonstrates the Von Neumann condition. Using the procedure outlined in Section \ref{sec:overset}, we obtain a very accurate solution for each of the wedge angles and obtain the transition criterion clearly. Note that the initial conditions mentioned above are used to obtain the coarse grid solution first and the coarse grid solution is used as an initial condition to get the fine grid solution using $\mathbf{P}^{4}$ based DGM. To show that the transition criterion has been captured accurately, we show the solution obtained for Mach 3.0 flow over a wedge with wedge angles $\theta_{w}=21.45^{\circ}$ and $\theta_{w}=21.46^{\circ}$ with the case 1 initial conditions in Figures \ref{fig:M3p0Angle21p45DFinal} and \ref{fig:M3p0Angle21p46DFinal} respectively.
\\

\noindent We also show the table containing the Mach stem height obtained using this procedure for various wedge angles for both cases of initial conditions in Table \ref{tab:1}. From Table \ref{tab:1}, we can see that the transition criterion occurs between wedge angles $21.45^{\circ}$ and $21.46^{\circ}$ and the Von Neumann transition criterion occurs between wedge angles $19.6^{\circ}$ and $19.7^{\circ}$. From the three-shock theory, the transition criterion is at $\theta_{w}=21.456^{\circ}$ and the Von Neumann condition occurs at $\theta_{w}=19.656^{\circ}$. This closely agrees with what we have obtained and validates our procedure.
\\
\begin{table}[h!]
    \caption{Non-dimensionalised height of the Mach stem ($H_{m}/H$) as a function of wedge angle ($\theta_w$) for $M$=3.0 and $w/H=1.0$ using two sets of initial conditions:- Case 1: Impulsive start, Case 2: Converged solution for $\theta_{w}=24^{\circ}$}
    \label{tab:1}
    \begin{center}
        \begin{tabular}{p{1.5cm}p{3.4cm}p{4.4cm}}
            \hline
         \multirow{2}{*}{} $\theta_{w}$ & $H_{m}/H$ for case 1 & $H_{m}/H$ for case 2 \\
          &  (impulsive start) & (converged solution for $\theta_{w}=24^{\circ}$) \\
            \hline
            $24^{\circ}$ & 0.274 & - \\
            \hline
            $23.5^{\circ}$ & 0.233 & 0.233 \\
            \hline
            $23^{\circ}$ & 0.188 & 0.188 \\
            \hline
            $22.5^{\circ}$ & 0.153 & 0.153 \\
            \hline
            $22^{\circ}$ & 0.122 & 0.122 \\
            \hline
            $21.46^{\circ}$ & 0.078 & 0.078 \\
            \hline
            $21.45^{\circ}$ & RR & 0.066 \\
            \hline
            $21^{\circ}$ & RR & 0.041 \\
            \hline
            $20.5^{\circ}$ & RR & 0.028 \\
            \hline
            $20^{\circ}$ & RR & 0.009 \\
            \hline
            $19.7^{\circ}$ & RR & 0.002 \\
            \hline
            $19.6^{\circ}$ & RR & RR \\
            \hline
            $19.5^{\circ}$ & RR & RR \\
            \hline
        \end{tabular}
    \end{center}
\end{table}

\noindent \textbf{2) Mach number 4.0:} We repeat the same procedure for various wedge angles between $\theta_{w}=20.5^{\circ}$ and $\theta_{w}=27^{\circ}$ near the transition criterion (transition from Regular reflection to Mach reflection) for Mach number 4.0 and $w/H=1.0$ in the computational domain shown in Figure \ref{fig:CompDomain}. Again, we consider two cases to demonstrate the hysteresis phenomenon. We solve the equations using an impulsive start as the initial condition for the first case, and the converged solution for $\theta_{w}=27^{\circ}$ as the second case. We show the table containing the Mach stem height obtained using this procedure for various wedge angles for both cases of initial conditions in Table \ref{tab:2}. From Table \ref{tab:2}, we can see that the transition criterion occurs between wedge angles $25.6^{\circ}$ and $25.7^{\circ}$ and the Von Neumann transition criterion occurs between wedge angles $20.9^{\circ}$ and $20.8^{\circ}$. From the three-shock theory, the transition criterion for Mach number $4.0$ is at $\theta_{w}=25.61^{\circ}$ and the Von Neumann condition occurs at $\theta_{w}=20.86^{\circ}$. This closely agrees with what we have obtained and validates our procedure.

\begin{table}[h!]
    \caption{Non-dimensionalised height of the Mach stem ($H_{m}/H$) as a function of wedge angle ($\theta_w$) for $M$=4.0 and $w/H=1.0$ using two sets of initial conditions:- Case 1: Impulsive start, Case 2: Converged solution for $\theta_{w}=27^{\circ}$}
    \label{tab:2}
    \begin{center}
        \begin{tabular}{p{1.5cm}p{3.4cm}p{4.4cm}}
            \hline
         \multirow{2}{*}{} $\theta_{w}$ & $H_{m}/H$ for case 1 & $H_{m}/H$ for case 2 \\
          &  (impulsive start) & (converged solution for $\theta_{w}=27^{\circ}$) \\
            \hline
            $27^{\circ}$ & 0.344 & - \\
            \hline
            $26.5^{\circ}$ & 0.313 & 0.313 \\
            \hline
            $26^{\circ}$ & 0.285 & 0.285 \\
            \hline
            $25.7^{\circ}$ & 0.245 & 0.245 \\
            \hline
            $25.6^{\circ}$ & RR & 0.224 \\
            \hline
            $25.5^{\circ}$ & RR & 0.203 \\
            \hline
            $25^{\circ}$ & RR & 0.183 \\
            \hline
            $24.5^{\circ}$ & RR & 0.144 \\
            \hline
            $24^{\circ}$ & RR & 0.115 \\
            \hline
            $23^{\circ}$ & RR & 0.093 \\
            \hline
            $22^{\circ}$ & RR & 0.045 \\
            \hline
            $21^{\circ}$ & RR & 0.005 \\
            \hline
            $20.9^{\circ}$ & RR & 0.002 \\
            \hline
            $20.8^{\circ}$ & RR & RR \\
            \hline
            $20.5^{\circ}$ & RR & RR \\
            \hline
        \end{tabular}
    \end{center}
\end{table}

\section{Conclusion}\label{sec:conc}
We have demonstrated the transition between regular reflection (RR) and Mach reflection (MR) of steady shock waves in inviscid flows using discontinuous Galerkin method (DGM) along with overset grids using an accurate shock capturing procedure \cite{srspkmr1}. We have identified the detachment criterion and the Von Neumann condition accurately and demonstrated the hysteresis which occurs in the transition. We have also calculated the Mach stem height for various wedge angles for Mach numbers 3.0 and 4.0. As future work, this shock capturing procedure can be used to further study the different flow phenomena that occur when the inflow Mach number is small ($<2.0$).


\bibliographystyle{ieeetr}
\bibliography{sn-bibliography}

\end{document}